\documentclass[10pt,a4paper,oneside]{amsart}
\usepackage{amsfonts, amscd, amsmath, amssymb}
\swapnumbers
\newtheorem{defn}[subsection]{Definition}

\newtheorem{lem}[subsubsection]{Lemma}
\newtheorem{cor}[subsubsection]{Corollary}
\newtheorem{rem}[subsubsection]{Remark}
\newtheorem{thm}[subsection]{Theorem}

\makeatletter
\@addtoreset{equation}{section}
\makeatother

\newcommand\RRR{\mathbb R}
\newcommand\va{a}
\newcommand\vb{b}
\newcommand\vn{n}
\newcommand\vm{m}
\newcommand\vi{i}
\newcommand\vj{j}
\newcommand\vk{k}

\newcommand\vt{t}
\newcommand\vx{x}

\newcommand\bnbh{V}

\newcommand\flow{\Phi}
\newcommand\aflow{\Phi}
\newcommand\bflow{\Psi}

\newcommand\Fld{F}
\newcommand\aFld{F}
\newcommand\bFld{G}

\newcommand\Fol{\mathcal F}
\newcommand\eps{\varepsilon}
\newcommand\Ieps{(-\eps,\eps)}

\newcommand\manif{M}

\newcommand\afunc{\alpha}
\newcommand\bfunc{\beta}

\newcommand\amap{f}
\newcommand\bmap{g}

\newcommand\difM{h}

\newcommand\der[2]{\frac{\partial #2}{\partial #1}}

\newcommand\dff[2]{#2.#1}
\newcommand\dffi[2]{\sum\limits_{\vi=1}^{\vm} #2^{\vi} \,\der{\vx_{\vi}}{#1}}

\newcommand\dfftotal[2]{D[#2;#1]} 

\newcommand\id{\mathrm{id}}

\newcommand\dAdx[1]{\frac{d#1}{d\vx}}
\newcommand\pdAdx[1]{\frac{\partial #1}{\partial\vx}}
\newcommand\pdAdt[1]{\frac{\partial #1}{\partial\vt}}

\newcommand\VA[1]{A_{#1}}
\newcommand\VF[1]{F_{#1}}
\newcommand\Va[1]{a^{#1}}
\newcommand\Vf[1]{f^{#1}}
\newcommand\spr[2]{\langle #1,#2\rangle}
\newcommand\sprFA[2]{\spr{\VF{#1}}{\VA{#2}}}

\newcommand\GD[1]{\mathbf{G}(#1)}
\newcommand\DFA[2]{D[#1;#2]}
\newcommand\prfa[2]{\Vf{#1}\Va{#2}}
\newcommand\sprFAi[3]{\sum\limits_{#3=1}^{\vm} \Vf{#3}_{#1} \Va{#3}_{#2}}

\newcommand\sprFAdn[2]{\sum\limits_{\vj=1}^{\vn} \Vf{#1}_{\vj} \Va{#2}_{\vj}}
\newcommand\sprFAup[2]{\sum\limits_{\vi=1}^{\vm} \Vf{\vi}_{#1} \Va{\vi}_{#2}}

\newcommand\mX{X}

\newcommand\mE{E}
\newcommand\mEm{\mE_{\vm}}
\newcommand\mEn{\mE_{\vn}}
\newcommand\mY{Y}

\newcommand\MA{\mathcal{A}}
\newcommand\MF{\mathcal{F}}

\newcommand\matrA{A}
\newcommand\matrB{B}

\newcommand\chpol[1]{P_{#1}}
\newcommand\coefX[1]{\mu_{#1}}

\newcommand\tran[1]{#1^{t}}

\newcommand\matrAe{\bar\matrA}
\newcommand\matrBe{\bar\matrB}

\newcommand\Shift{\phi}
\newcommand\smr{C^{\infty}(\manif,\RRR)}
\newcommand\smm{C^{\infty}(\manif,\manif)}
\newcommand\Fldset{D}

\newcommand\Mn{M(\vn)}
\newcommand\Mm{M(\vm)}
\newcommand\Mmn{M(\vm,\vn)}

\newcommand\Mnn{M(\vn,\vn)}
\newcommand\chpolmap{p}
\newcommand\difchpol{\delta}\begin{document}
\author{Sergey Maksymenko}
\title{Consecutive shifts along orbits of vector fields}

\address{
 Topology Department,
 Institute of Mathematics, NAS of Ukraine,
 Tereshchenkivska str. 3, 01601 Kyiv, Ukraine,
 e-mail:\texttt{maks@imath.kiev.ua},
 phone: (+380 44) 2345150,
 fax: (+380 44) 2352010 
}

\begin{abstract}
Let $M$ be a smooth ($C^{\infty}$) manifold, $F_1,...,F_n$ be vector fields on $M$ generating the corresponding flows 
$\Phi_1,...,\Phi_n$, and $\alpha_1,\ldots,\alpha_{n}:M\to \mathbb{R}$ smooth functions.
Define the following map $f:M\to M$ by  
$$f(x)= \Phi_n ( ...  ( \Phi_2 ( \Phi_1 (x,\alpha_1(x) ), \alpha_2(x) ),  ... , \alpha_n(x) ).$$ 
In this note we give a necessary and sufficient condition on vector fields $F_1,...,F_n$ and smooth functions  $\alpha_1,\ldots,\alpha_{n}$
for $f$ to be a local diffeomorphism.
\end{abstract}

\maketitle
\section{Introduction}
Let $\manif$ be a smooth ($C^{\infty}$) manifold and $\Fldset$ be an arbitrary set of vector fields on $\manif$.
Following~\cite{Sussmann:BAMS} say that two points $\vx$ and $\vx'$ are $\Fldset$-connected provided there exist vector fields $\Fld_1,\ldots,\Fld_{\vn} \in\Fldset$ and real numbers $\vt_1,\ldots\vt_{\vn}$ such that 
\begin{equation}\label{equ:Orbz}
 \vx' = \flow_{\vn}( \cdots  (\flow_{2}( \flow_{1}(\vx, \vt_{1}), \vt_{2}),  \cdots  , \vt_{\vn}),
\end{equation}
where $\flow_{\vi}$ is a (local) flow generated by $\Fld_{\vi}$.
This defines an equivalence relation and the corresponding equivalence classes are called {\em orbits}.
The partition of $\manif$ by such orbits is a foliation $\Fol$ with singularities on $\manif$, see~\cite{Sussmann:BAMS,Sussmann:TransAMS,Stefan:BAMS,Stefan:JLMS}.

In this note we consider mappings obtained by consecutive smooth shifts along orbits of vector fields belonging to $\Fldset$.

For $\Fld_1,\ldots,\Fld_{\vn} \in\Fldset$ and $\afunc_1,\ldots,\afunc_{\vn}\in\smr$ define the following map 
$$\Shift_{1, \ldots, \vn}(\afunc_1,\ldots,\afunc_{\vn}):\manif\to\manif$$
by
\begin{equation}\label{equ:shift_def}
\Shift_{1, \ldots, \vn}(\afunc_1,\ldots,\afunc_{\vn})(x)\ = \ \flow_{\vn}\bigl( \cdots  (\flow_{2}( \flow_{1}(\vx, \afunc_1(x)), \afunc_2(\vx)),  \cdots  \afunc_{\vn}(\vx)\bigr). 
\end{equation}
We shall say that this map is a {\em shift-map along $\flow_1,\ldots,\flow_{\vn}$ via functions $\afunc_1,\ldots,\afunc_{\vn}$}.
These functions will be called {\em shift-functions\/} for 
the mapping $\Shift_{1, \ldots, \vn}(\afunc_1,\ldots,\afunc_{\vn})$.

\begin{rem}\rm
Notice that if $\manif$ is non-compact, then in general a vector field $\Fld$ on $\manif$ does not generate a global flow.
However, if $\bnbh\subset\manif$ is open and has a compact closure, then $\Fld$ generates a local flow $\flow:\bnbh\times(-\eps,\eps)\to\manif$ for some $\eps>0$.
Thus the shift-mapping $\Shift_{1, \ldots, \vn}(\afunc_1,\ldots,\afunc_{\vn})$ can not be defined for arbitrary functions $\afunc_1,\ldots,\afunc_{\vn}$.
However, throughout this paper when speaking about a shift-map 
$\Shift_{1, \ldots, \vn}(\afunc_1,\ldots,\afunc_{\vn})$
we shall always assume that it is indeed well-defined on all of $\manif$.
\end{rem}

\begin{rem}\label{rem:sing_coord}\rm
In a certain sense, we may imagine that vector fields $\Fld_{\vj}$ define some ``singular coordinate system'' on $\manif$ so that the functions $\afunc_{\vi}$ are ``coordinate functions'' of the mapping  $\Shift_{1, \ldots, \vn}(\afunc_1,\ldots,\afunc_{\vn})$ in this ``coordinate system.''
It turns out, see Remark~\ref{rem:Fol_manif}, that if this ``coordinate system'' is ``right'' then $\afunc_{\vi}$ are just differences between ``right'' coordinate functions of $\amap$ and ``right'' coordinates.
\end{rem}

Evidently, $\Shift_{1, \ldots, \vn}(\afunc_1,\ldots,\afunc_{\vn})$ preserves each leaf of the foliation $\Fol$.  
Conversely, in some cases, e.g.\! when $\Fol$ is a non-singular foliation, each leaf-preserving diffeomorphism $\difM:\manif\to\manif$ can be represented (in general only locally!) as a smooth shift via some smooth functions, see Section~\ref{sect:repr_sm_shifts}.

Our main result is Theorem~\ref{th:f_loc_dif} that gives a necessary and sufficient condition on the functions $\afunc_1,\ldots,\afunc_{\vn}$ and vector fields $\Fld_1,\ldots,\Fld_{\vn}$ for the map 
$\Shift_{1, \ldots, \vn}(\afunc_1,\ldots,\afunc_{\vn})$ to be a local diffeomorphism at some point.
In fact, we obtain the expression for the Jacobian of $\amap$ in the (invariant) terms of derivatives of $\afunc_{\vi}$ along $\Fld_{\vj}$.
In the author's paper~\cite{Maks:Shifts} the similar result was obtained for the shift-mappings along the orbits of one flow.

The unexpected feature of the obtained condition is that it is invariant with respect to the simultaneous permutations of indices. 
Thus, if $\sigma\in\Sigma_{\vn}$ is a permutation of $1,\ldots,\vn$, then it turns out (see Lemma~\ref{lem:thm_inv_perm}) that 
the mapping $\Shift_{1,\ldots,\vn}(\afunc_1,\ldots,\afunc_{\vn})$ is a diffeomorphism of $\manif$ iff so is $$\Shift_{\sigma(1),\ldots,\sigma(\vn)}(\afunc_{\sigma(1)},\ldots,\afunc_{\sigma(\vn)}),$$
see also Corollary~\ref{cor:commut}.

\begin{rem}\em
Similarly to~\eqref{equ:shift_def} for every flow $\flow_{\vi}$ we can also define the {\em shift-map\/} 
$$\Shift_i:\smr\to\smm$$  by $\Shift_i(\afunc)(x)=\flow(x,\afunc(x))$, see~\cite{Maks:Shifts}.

Emphasize that in~\eqref{equ:shift_def} we take the values of each $\afunc_{\vi}$ {\em at the initial point $x$}.
Therefore, in general, the mapping $\Shift_{1, \ldots, \vk}(\afunc_1,\ldots,\afunc_{\vk})$ does not coincide with the composition 
$\Shift_{\vk}(\afunc_{\vk}) \circ \cdots \circ \Shift_{2}(\afunc_{2}) \circ \Shift_{1}(\afunc_{1})$, e.g. for $\vn=2$ we have:
$$
\Shift_{1,2}(\afunc_1,\afunc_{2})(\vx)=
\flow_{2}\bigl(\flow_{1}\bigl(\vx,\afunc_{1}(\vx)\bigr),\afunc_{2}(\vx)\bigr)
$$
while
$$
\Shift_{2}(\afunc_{2}) \circ \Shift_{1}(\afunc_{1})(\vx)=
\flow_{2}\bigl(
 \flow_{1}\bigl(\vx,\afunc_{1}(\vx)\bigr),\
 \afunc_{2}\circ\flow_{1}\bigl(\vx,\afunc_{1}(\vx)\bigr)\bigr).
$$

In particular, it easily follows that if two adjacent vector fields coincide, say $\flow_{\vi}=\flow_{\vi+1}$, then 
\begin{equation}\label{equ:adj_ind_equal}
\Shift_{1,\ldots, \vi, \vi, \ldots, \vn}(\afunc_1,\ldots, \afunc'_{\vi}, \afunc''_{\vi},\ldots,\afunc_{\vn})
=
\Shift_{1,\ldots, \vi, \ldots, \vn}(\afunc_1,\ldots, \afunc'_{\vi}+\afunc''_{\vi}, \ldots,\afunc_{\vn}).
\end{equation}
\end{rem}

\subsection{Structure of the paper.}
In Section~\ref{sect:ext_Gramm_det} for each pair of $\vn$-tuples of vectors in $\RRR^{\vm}$ we define two square matrices $\mX$ and $\mY$ of dimensions $\vm\times\vm$ and $\vn\times\vn$ respectively such that $|\mEm+\mX|=|\mEn+\mY|$, where $\mEm$ and $\mEn$ are unit matrices, see Definition~\ref{def:symbol_D}.

Further, in Section~\ref{sect:Shift_Are_Diffs} we show that the latter equality gives us two expressions of the Jacobi determinant of a shift-mapping at its fixed point (Theorem~\ref{th:f_loc_dif}).
In fact, $|\mE+\mX|$ is the usual Jacobian, while $|\mE+\mY|$ is the Jacobian with respect to ``singular coordinates'', it depends only on the derivatives of shift functions $\afunc_{\vi}$ along vector fields $\Fld_{\vj}$, see Remark~\ref{rem:sing_coord}.

Finally in Section~\ref{sect:repr_sm_shifts} we give examples of foliations generated by certain vector fields $\Fld_{1},\ldots,\Fld_{\vn}$ such that every leaf-preserving mapping can be represented as a shift along the orbits of these vector fields via some smooth functions.

\subsection{Acknowledgements}
I would like to thank V.~V.~Sharko, V.~V.~Sergeichuk, D.~Bolotov, and E.~Polulyah for valuable conversation.

\section{Characteristic polynomials of products of two matrices}\label{sect:ext_Gramm_det}
 
We shall designate by $\Mmn$ the space of all $\vm\times\vn$-matrices (e.g. matrices with $\vm$ rows and $\vn$ columns).
If $\vm=\vn$, then $\Mnn$ will be denoted by $\Mn$.
For each $\mX\in\Mn$ let 
\begin{equation}\label{equ:chpolX}
\chpol{\mX}(\lambda) = |\mX-\lambda\mEn| = (-\lambda)^{\vm} + 
(-\lambda)^{\vm-1}\coefX{1}  +  (-\lambda)^{\vm-2}\coefX{2} + \ldots + \coefX{\vm}
\end{equation}
be the characteristic polynomial of $\mX$.
Denote by $\mEn$ the identity matrix of dimension $\vm$, and
by $0_{\vn,\vm}$ the zero $(\vn\times\vm)$-matrix.

\begin{thm}
Let $\matrA,\matrB\in\Mmn$, thus the matrix $\matrA\tran{\matrB}$ is of dimension $\vm\times\vm$ and
$\tran{\matrA}\matrB$ is of dimension $\vn\times\vn$.
Then
$$
\chpol{\matrA\tran{\matrB}}(\lambda)=(-\lambda)^{\vm-\vn}\chpol{\tran{\matrA}\matrB}(\lambda).
$$
\end{thm}
\proof
Consider 3 cases.

{\bf Case 1.}~If $\vm=\vn$, then theorem claims that 
$\chpol{\matrA\tran{\matrB}}=\chpol{\tran{\matrA}\matrB}$.
This is implied by the following lemma:
\begin{lem}\label{lem:pAB_pBA}
Let $\matrA,\matrB\in\Mn$. Then
$\chpol{\matrA\matrB}(\lambda)\equiv\chpol{\matrB\matrA}(\lambda)$.
\end{lem}
Indeed,
$$
\chpol{\matrA\tran{\matrB}}=
\chpol{\tran{(\matrA\tran{\matrB})}}=
\chpol{\matrB\tran{\matrA}}
\stackrel{\text{Lemma~\ref{lem:pAB_pBA}}}{=\!=\!=\!=\!=\!=\!=\!=\!=}
\chpol{\tran{\matrA}\matrB}.
$$

\proof[Proof of Lemma~\ref{lem:pAB_pBA}.]
This result is known.
But for the completeness and for the convenience of the reader we present a short topological proof.

Suppose that one of the matrices, say $\matrA$ is non-degenerate.
Then the identity $\matrA(\matrB\matrA)=(\matrA\matrB)\matrA$ implies that the matrices $\matrA\matrB$ and $\matrB\matrA$ are conjugate:
\begin{equation}\label{equ:BA_cong_AB}
\matrB\matrA = \matrA^{-1}(\matrA\matrB)\matrA,
\end{equation}
whence $\chpol{\matrA\matrB}=\chpol{\matrB\matrA}$.

If both $\matrA$ and $\matrB$ are degenerated, then it is possible that $\matrA\matrB$ and $\matrB\matrA$ are not conjugate.
For example, if 
$$
\matrA=\left(\begin{matrix} 1 & 0 \\ 0 & 0 \end{matrix} \right) \qquad\text{and}\qquad
\matrB=\left(\begin{matrix} 0 & 1 \\ 0 & 0 \end{matrix} \right),
$$
then 
$\matrA\matrB = \left(\begin{matrix} 0 & 1 \\ 0 & 0 \end{matrix} \right)$, while
$\matrB\matrA = \left(\begin{matrix} 0 & 0 \\ 0 & 0 \end{matrix} \right)$.
However, $\matrA\matrB$ and $\matrB\matrA$ always have the same characteristic polynomials due to the following arguments.

Define the following mappings 
$$
\begin{array}{lcl}
\chpolmap:\Mn\to\RRR^{\vn} & \qquad\qquad &
\chpolmap(\mX)=(\coefX{1},\coefX{2},\ldots,\coefX{\vn}), \\ [2mm]
\difchpol:\Mn\times\Mn\to\RRR^{\vn} & \qquad\qquad &
\difchpol(\matrA,\matrB) = \chpolmap(\matrA\matrB)-\chpolmap(\matrB\matrA),
\end{array}
$$
where $\coefX{\vi}$ are given by~\eqref{equ:chpolX}.
By~\eqref{equ:BA_cong_AB} $\difchpol(\matrA,\matrB)=0$ provided at least one of the matrices either $\matrA$ or $\matrB$ is non-degenerate.
Since $\difchpol$ is continuous and the subset of $\Mn\times\Mn$ consisting of pairs $(\matrA,\matrB)$ in which both $\matrA$ and $\matrB$ are degenerated is nowhere dense, it follows that $\difchpol\equiv0$ on $\Mn\times\Mn$, i.e. $\chpol{\matrA\matrB}=\chpol{\matrB\matrA}$ for all $\matrA,\matrB\in\Mn$.
\qed

{\bf Case 2.}~Suppose that $\vm>\vn$.
Let us add to $\matrA$ and $\matrB$ $\vm-\vn$ zero columns and denote the obtained $\vm\times\vm$-matrices by $\matrAe$ and $\matrBe$.
Then it is easy to see that 
$$
\matrAe\tran{\matrBe} = \matrA\tran{\matrB}
\qquad\text{and}\qquad
\tran{\matrAe}\matrBe=
\left(
\begin{matrix}
\matrAe\tran{\matrBe} & 0_{\vn,\vm-\vn} \\
0_{\vm-\vn,\vn} & 0_{\vm-\vn,\vm-\vn}
\end{matrix}
\right).
$$
Therefore
$$
\chpol{\matrA\tran{\matrB}}(\lambda) \ =  \ \chpol{\matrAe\tran{\matrBe}}(\lambda)
\ \stackrel{\text{Case 1}}{=\!=\!=\!=\!=} \
\chpol{\tran{\matrAe}\matrBe}(\lambda) \ = \ (-\lambda)^{\vm-\vn}\chpol{\tran{\matrA}\matrB}(\lambda).
$$

{\bf Case 3.}~ The case $\vm<\vn$ reduces to the case 2 by transposing $\matrA$ and $\matrB$.
Theorem is proved.
\endproof

As a corollary we obtain the following identity which will play the crucial role:
\begin{equation}\label{equ:E_AB__E_BA}
|\mEm+\matrA\tran{\matrB}|=\chpol{\matrA\tran{\matrB}}(-1)=\chpol{\tran{\matrA}\matrB}(-1)=|\mEn+\tran{\matrA}\matrB|.
\end{equation}

\subsection{Symbol $\DFA{-}{-}$.}
For two vectors 
$$
\VF{}=
\left(
\begin{array}{c}
\Vf{1} \\ \Vf{2} \\ \cdots \\ \Vf{\vm}
\end{array}
\right)
  \qquad \text{and} \qquad 
\VA{}=
\left(
\begin{array}{c}
\Va{1} \\ \Va{2} \\ \cdots \\ \Va{\vm}
\end{array}
\right)
$$
we can define the following two products:
$$ 
\VF{}\tran{\VA{}} = 
\left(
\begin{array}{c}
\Vf{1} \\ \Vf{2} \\ \cdots \\ \Vf{\vm}
\end{array}
\right)
\cdot
\left(\Va{1},\ldots,\Va{\vm}\right)
=
\left(
\begin{array}{cccc}
\prfa{1}{1}   & \prfa{1}{2}   & \cdots & \prfa{1}{\vm}   \\
\prfa{2}{1}   & \prfa{2}{2}   & \cdots & \prfa{2}{\vm}   \\
\cdots         & \cdots         & \cdots & \cdots           \\
\prfa{\vm}{1} & \prfa{\vm}{2} & \cdots & \prfa{\vm}{\vm} 
\end{array}
\right) \ \in \ \Mm
$$
and
$$\tran{\VF{}} \VA{} = \spr{\VF{}}{\VA{}} = \sprFAi{}{}{\vi} \ \in \ \RRR. $$ 

Fix two systems of vectors $\VF{1},\ldots,\VF{\vn}$ and $\VA{1},\ldots,\VA{\vn}$ in $\RRR^{\vm}$ and define the following two matrices:

\begin{equation}\label{equ:def_X_non_coord}
\mX = \VF{1} \tran{\VA{1}} + \ldots + \VF{\vn} \tran{\VA{\vn}}
\end{equation}
and 
\begin{equation}\label{equ:def_G_non_coord}
\mY =
\left(
\begin{array}{cccc}
\sprFA{1}{1}   & \sprFA{2}{1}   & \cdots & \sprFA{\vn}{1}   \\
\sprFA{1}{2}   & \sprFA{2}{2}   & \cdots & \sprFA{\vn}{2}   \\
\cdots         & \cdots         & \cdots & \cdots           \\
\sprFA{1}{\vn} & \sprFA{2}{\vn} & \cdots & \sprFA{\vn}{\vn} 
\end{array}    
\right).
\end{equation}
Notice that if $\VF{\vi}=\VA{\vi}$ for all $\vi=1,\ldots,\vn$, then $|\mY|$ is the Gramm determinant $\GD{\VF{1},\ldots,\VF{\vn}}$.

Let $\VF{\vi}=(\Vf{1}_{\vi},\ldots,\Vf{\vm}_{\vi})$ and 
$\VA{\vi}=(\Va{1}_{\vi},\ldots,\Va{\vm}_{\vi})$ be the coordinates of these vectors and 
$$ 
\MF = 
 \left(
\begin{array}{cccc}
\Vf{1}_{1}   & \Vf{1}_{2}   & \cdots & \Vf{1}_{\vn}   \\ [2mm]
\Vf{2}_{1}   & \Vf{2}_{2}   & \cdots & \Vf{2}_{\vn}   \\ [2mm]
\cdots           & \cdots           & \cdots & \cdots           \\
\Vf{\vm}_{1} & \Vf{\vm}_{2} & \cdots & \Vf{\vm}_{\vn} 
\end{array}    
\right),
\qquad
\MA = 
 \left(
\begin{array}{cccc}
\Va{1}_{1}   & \Va{1}_{2}   & \cdots & \Va{1}_{\vn}   \\ [2mm]
\Va{2}_{1}   & \Va{2}_{2}   & \cdots & \Va{2}_{\vn}   \\ [2mm]
\cdots           & \cdots           & \cdots & \cdots           \\
\Va{\vm}_{1} & \Va{\vm}_{2} & \cdots & \Va{\vm}_{\vn} 
\end{array}    
\right)
$$
be the matrices whose columns consist of coordinates of $\VF{\vi}$ and $\VA{\vi}$.
A simple calculation shows that 
\begin{equation}\label{equ:def_X}
\mX =  \MF \MA^{t} = 
 \left(
\begin{array}{cccc}
\sprFAdn{1}{1}   & \sprFAdn{1}{2}   & \cdots & \sprFAdn{1}{\vm}   \\ [2mm]
\sprFAdn{2}{1}   & \sprFAdn{2}{2}   & \cdots & \sprFAdn{2}{\vm}   \\ [2mm]
\cdots           & \cdots           & \cdots & \cdots           \\
\sprFAdn{\vm}{1} & \sprFAdn{\vm}{2} & \cdots & \sprFAdn{\vm}{\vm} 
\end{array}    
\right) \ \in \ \Mm
\end{equation}
and 
\begin{equation}\label{equ:def_G}
\mY = \MF^{t} \MA =
 \left(
\begin{array}{cccc}
\sprFAup{1}{1}   & \sprFAup{1}{2}   & \cdots & \sprFAup{1}{\vn}   \\ [2mm]
\sprFAup{2}{1}   & \sprFAup{2}{2}   & \cdots & \sprFAup{2}{\vn}   \\ [2mm]
\cdots           & \cdots           & \cdots & \cdots           \\
\sprFAup{\vn}{1} & \sprFAup{\vn}{2} & \cdots & \sprFAup{\vn}{\vn} 
\end{array}    
\right) \ \in \ \Mn.
\end{equation}

Then it follows from~\eqref{equ:E_AB__E_BA} for $\matrA=\MF$ and $\matrB=\MA$ that 
$$
|\mEm+\mX| = |\mEn+\mY|.
$$

\begin{defn}\label{def:symbol_D}
For a pair of $\vn$-tuples of vectors $$\VF{1},\ldots,\VF{\vn}\qquad\text{and}\qquad\VA{1},\ldots,\VA{\vn}$$ in $\RRR^{\vm}$ let us define the following symbol $\DFA{\VF{1},\ldots,\VF{\vn}}{\VA{1},\ldots,\VA{\vn}}$ by
\begin{equation}\label{equ:DFA}
\DFA{\VF{1},\ldots,\VF{\vn}}{\VA{1},\ldots,\VA{\vn}}:=|\mEm+\mX| = |\mEn+\mY|,
\end{equation}
where $\mX$ and $\mY$ are given by formulas~\eqref{equ:def_X_non_coord} and~\eqref{equ:def_G_non_coord} respectively.
\end{defn}

\subsection{Functions and vector fields.}
Let $\manif$ be a smooth manifold.
For a smooth vector field $\aFld$ on $\manif$ and a smooth function $\afunc:\manif\to\RRR$ the derivative of $\afunc$ along $\aFld$ will be denoted by $\dff{\afunc}{\aFld}$:
$$
\dff{\afunc}{\aFld} = \dffi{\afunc}{\aFld} = \spr{\aFld}{\nabla\afunc}.
$$
 
Let $\aFld_1,\ldots,\aFld_{\vn}$ be vector fields and $\afunc_1,\ldots,\afunc_{\vn}$ be smooth functions on $\manif$.
Then at each point of $\manif$ we have two systems of vectors 
$\aFld_1,\ldots,\aFld_{\vn}$ and $\nabla\afunc_1,\ldots,\nabla\afunc_{\vn}$.
Hence the following symbol 
$$\DFA{\VF{1},\ldots,\VF{\vn}}{\nabla\afunc_1,\ldots,\nabla\afunc_{\vn}}$$
is well defined by formula~\eqref{equ:DFA}.
For simplicity we shall denote it by \begin{equation}\label{equ:dfftotal}
\dfftotal{\afunc_1,\ldots,\afunc_{\vn}}{\aFld_1,\ldots,\aFld_{\vn}} := 
\DFA{\VF{1},\ldots,\VF{\vn}}{\nabla\afunc_1,\ldots,\nabla\afunc_{\vn}}.
\end{equation}

Formula~\eqref{equ:DFA} gives two ways for calculation of $\DFA{\VF{1},\ldots,\VF{\vn}}{\VA{1},\ldots,\VA{\vn}}$.
The first expressions $|\mEm+\mX|$ depend on local coordinates.
However the second one $|\mEn+\mY|$ includes only derivatives of $\afunc_{\vi}$ along $\Fld_{\vj}$, and thus is invariant with respect to the local coordinates, whence so is $\DFA{\VF{1},\ldots,\VF{\vn}}{\VA{1},\ldots,\VA{\vn}}$.

In particular, if $\aFld$ and $\bFld$ are smooth vector fields and $\afunc$ and $\bfunc$ are smooth functions on $\manif$ then 
\begin{multline}\label{equ:1_flow} 
\dfftotal{\afunc}{\aFld} = 1 + \dff{\afunc}{\aFld},\hfill
\end{multline}
\begin{multline}\label{equ:2_flow}  
\dfftotal{\afunc,\bfunc}{\aFld,\bFld} = 
\left|
\begin{array}{cc}
 1+\dff{\afunc}{\aFld} & \dff{\bfunc}{\aFld} \\
 \dff{\afunc}{\bFld} & 1+\dff{\bfunc}{\bFld}
\end{array}
\right|=1 + \dff{\afunc}{\aFld}  + \dff{\bfunc}{\bFld}  + 
\left|
\begin{array}{cc}
 \dff{\afunc}{\aFld} & \dff{\bfunc}{\aFld} \\
 \dff{\afunc}{\bFld} & \dff{\bfunc}{\bFld}
\end{array}
\right|,\hfill
\end{multline} 
and so on.

Notice that~\eqref{equ:1_flow} coincides with formula~(14) of~\cite{Maks:Shifts} for the Jacobi determinant of a shift-mapping along the trajectories of $\Fld$ via function $\afunc$.

\begin{lem}\label{lem:thm_inv_perm}
Let $\sigma\in\Sigma_{\vn}$ be a permutation of indices $1,\ldots,\vn$. Then 
$$
\dfftotal{\afunc_1,\ldots,\afunc_{\vn}}{\aFld_1,\ldots,\aFld_{\vn}}
=
\dfftotal{\afunc_{\sigma(1)},\ldots,\afunc_{\sigma(\vn)}}{\aFld_{\sigma(1)},\ldots,\aFld_{\sigma(\vn)}},
$$
\end{lem}
\proof
By formula~\eqref{equ:def_X_non_coord} for the matrix $\mX$ the simultaneous permutation of vector fields and functions does not change the matrix $\mE+\mX$:
$$
\mEm + \mX = \mEm + \sum\limits_{\vi=1}^{\vn} \VF{\vi}\cdot \tran{\nabla\afunc_{\vi}} = 
\mEm + \sum\limits_{\vi=1}^{\vn} \VF{\sigma(\vi)} \cdot \tran{\nabla\afunc_{\sigma(\vi)}}. 
$$

This lemma can also be established using matrix $\mY$.
It suffices to consider the case when $\sigma$ is a transposition $(\vi\vj)$. 
Then it yields a simultaneous exchanging of $\vi$-th and $\vj$-th columns and $\vi$-th and $\vj$-th rows. 
Such a procedure does not change the determinant $|\mEn+\mY|$. 
\endproof

\section{Shifts that are diffeomorphisms}\label{sect:Shift_Are_Diffs}
Let $\Fld_1,\ldots,\Fld_{\vn}$ be vector fields on $\RRR^{\vm}$, and for each $\vi=1,\ldots,\vn$ 
$$\flow_{\vi}:\bnbh\times(-\eps,\eps)\to\RRR^{\vm}$$ 
a local flow generated by $\Fld_{\vi}$ on some neighborhood $\bnbh$ of $0\in\RRR^{\vm}$.
Let $$\afunc_1,\ldots,\afunc_{\vn}:\RRR^{\vm}\to\RRR$$ be smooth functions such that the shift mapping
$\Shift_{1,\ldots,\vn}(\afunc_1,\ldots,\afunc_{\vn})$ is well-defined.
For simplicity, we shall denote this map by $\Shift(\afunc_1,\ldots,\afunc_{\vn})$.

In this section we give a necessary and sufficient condition on functions $\afunc_1,\ldots,\afunc_{\vn}$ for the corresponding shift-map $\Shift(\afunc_1,\ldots,\afunc_{\vn})$ to be a diffeomorphism.

Recall, see~\eqref{equ:dfftotal}, that we have introduced the following symbol $$\dfftotal{\afunc_1,\ldots,\afunc_{\vn}}{\aFld_1,\ldots,\aFld_{\vn}}$$
being a certain expression that included {\em only the derivatives of each $\afunc_{\vi}$ along each $\aFld_{\vj}$}.
In fact, it will turn out that if $\afunc_{\vi}(0)=0$ for all $\vi$, and thus $\amap(0)=0$, then this symbol is just the Jacobian of $\amap$ at $0$.
\begin{thm}\label{th:f_loc_dif}
The mapping $\amap=\Shift(\afunc_1,\ldots,\afunc_{\vn})$, defined by~\eqref{equ:shift_def}, is a local diffeomorphism at $0\in\RRR^{\vm}$ if and only if
\begin{equation}\label{equ:f_dif_cond}
\dfftotal{\afunc_1,\ldots,\afunc_{\vn}}{\aFld_1,\ldots,\aFld_{\vn}}\not=0. 
\end{equation}
Moreover, $\amap$ preserves orientation of $\RRR^{\vm}$ iff $\dfftotal{\afunc_1,\ldots,\afunc_{\vn}}{\aFld_1,\ldots,\aFld_{\vn}}>0$.
\end{thm}

From Lemma~\ref{lem:thm_inv_perm} we obtain the following corollary:
\begin{cor}\label{cor:thm_inv_perm}
Let $\sigma\in\Sigma_{\vn}$ be a permutation of indices $\{1,\ldots,\vn\}$.
Then the following shift-map $$\Shift_{\sigma(1),\ldots,\sigma(\vn)}(\afunc_{\sigma(1)},\ldots,\afunc_{\sigma(\vn)})$$
along $\flow_{\sigma(1)},\ldots,\flow_{\sigma(\vn)}$ is a preserving (reversing) orientation  diffeomorphism iff so is $\Shift_{1,\ldots,\vn}(\afunc_{1},\ldots,\afunc_{\vn})$
\end{cor}
Thus in order to establish that a shift mapping $\Shift_{1,\ldots,\vn}(\afunc_{1},\ldots,\afunc_{\vn})$ is a diffeomorphism we may replace it by another shift mapping simultaneously transposing corresponding flows and functions.

As a simple application we prove the following statement:
\begin{cor}\label{cor:commut}
For {\em arbitrary\/} vector fields $\Fld_1,\Fld_2$ and {\em arbitrary\/} smooth functions $\afunc_1,\afunc_2$ the ``commutator'' $\amap=\Shift_{1,2,1,2}(\afunc_1,\afunc_2,-\afunc_1,-\afunc_2)$ is {\em always\/} a (local) diffeomorphism, even if these vector fields do not commute.
\end{cor}
\proof
By Corollary~\ref{cor:thm_inv_perm}, this mapping is a local diffeomorphism iff the following map
$$
\Shift_{1,1,2,2}(\afunc_1,-\afunc_1,\afunc_2,-\afunc_2)
\ \stackrel{\eqref{equ:adj_ind_equal}}{\equiv\!\equiv\!\equiv} \
\Shift_{1,2}(0,0) \ \equiv \ \id
$$
is. But the identity mapping is of course a diffemorphism, whence so is $\amap$.
\endproof

\subsection{Proof of Theorem~\ref{th:f_loc_dif}.}
We may suppose that each $\afunc_{\vi}(0)=0$.
Otherwise, set $C_{\vi}=\afunc_{\vi}(0)$.
Regarding each $C_{\vi}:\bnbh\to\RRR$ as a constant function, we see that the smooth shift $\bmap=\Shift(C_{1},C_{2},\ldots,C_{\vn})$ via them is always a diffeomorphism.
Hence, $\amap$ is a diffeomorphism at $0\in\RRR^{\vm}$ iff $\bmap^{-1}\circ\amap$ is.
But it is easy to see that 
$$
\bmap^{-1}\circ\amap=\Shift(\afunc_{1}-C_{1},\afunc_{2}-C_{2},\ldots,\afunc_{\vn}-C_{\vn})
$$
is a smooth shift along $\Fol$ via the functions $\afunc_{\vi}-C_{\vi}$ vanishing at $0$.
Moreover, 
$$
\dfftotal{\afunc_1-C_1,\ldots,\afunc_{\vn}-C_{\vn}}{\aFld_1,\ldots,\aFld_{\vn}}
\ = \
\dfftotal{\afunc_1,\ldots,\afunc_{\vn}}{\aFld_1,\ldots,\aFld_{\vn}},
$$
since 
$$
\dff{(\afunc_{\vi}-C_{\vi})}{\aFld_{\vj}}=\dff{\afunc_{\vi}}{\aFld_{\vj}},
\qquad \vi,\vj=1,\ldots,\vn.
$$

Thus we can replace $\amap$ with $\bmap^{-1}\circ\amap$ and assume that $\afunc_{\vi}(0)=0$.
Then $\amap(0)=0$ and all we need is to calculate the Jacobian $J(\amap,0)=|\dAdx{\amap}(0)|$ of $\amap$ at $0$.
The following lemma implies our theorem.

\begin{lem}\label{lm:Jf0_E_X}
Let $\mEm$ be a unit $\vm\times\vm$-matrix.
If $\afunc_{\vi}(0)=0$ for all $\vi=1,\ldots,\vn$, then
$$
J(\amap,0) \ = \ 
|\mEm +  \Fld_1\cdot\tran{\nabla\afunc_1} + \cdots + \Fld_{\vn}\cdot\tran{\nabla\afunc_{\vn}}|
\ \stackrel{\eqref{equ:DFA}}{=\!=\!=} \
\dfftotal{\afunc_1,\ldots,\afunc_{\vn}}{\Fld_1,\ldots,\Fld_{\vn}}.
$$ 
\end{lem}
\proof
We shall assume that $\vn=2$.
The general case is quite analogous.
Notice that for $\vn=1$ this lemma coincides with Lemma~20 of~\cite{Maks:Shifts}, and that this case also follows from the case $\vn=2$ for $\Fld_{2}\equiv0$.

For the convenience, let us change the notations.
Thus suppose we have two vector fields $\aFld,\bFld$ on $\RRR^{\vm}$ generating two local flows 
$$
 \aflow=(\aflow^{1},\ldots,\aflow^{\va}), \ \bflow=(\bflow^{1},\ldots,\bflow^{\vb}) \ : \ \bnbh\times\Ieps\to\RRR^{\vm},
$$
and two smooth functions $\afunc,\bfunc:\bnbh\to\Ieps$ such that $\afunc(0)=\bfunc(0)=0$.
Then
$$ 
\amap(\vx) = \Shift(\afunc,\bfunc)(x)=
\bflow( \aflow(x,\afunc(x)), \bfunc(x) ).
$$

Differentiating coordinate functions of $\amap$ by $\vx_1,\ldots,\vx_{\vn}$ we get:

$$
\dAdx{\amap} =
\pdAdx{\bflow} \cdot \left(\pdAdx{\aflow} + \pdAdt{\aflow}\cdot \nabla\afunc \right) + 
\pdAdt{\bflow} \cdot\nabla\bfunc,
$$
where 
$$
\pdAdx{\aflow}=\left( \der{\vx^{j}}{\aflow^{i}} \right) 
\qquad \text{and} \qquad
\pdAdx{\bflow}=\left( \der{\vx^{j}}{\bflow^{i}} \right)
$$
are $\vm\times\vm$-matrices,
$$
\pdAdt{\aflow}=\left( \der{\vt}{\aflow^{i}} \right)
\qquad \text{and} \qquad
\pdAdt{\bflow}=\left( \der{\vt}{\bflow^{i}} \right)
$$
are $\vm$-vectors, and the expressions 
$$
\pdAdt{\aflow} \cdot \nabla\afunc = \left(\der{\vt}{\aflow^{i}}\right) \cdot \tran{\nabla\afunc} \qquad 
\text{and} \qquad
\pdAdt{\bflow} \cdot \nabla\bfunc = \left(\der{\vt}{\bflow^{i}}\right) \cdot \tran{\nabla\bfunc},
$$
where $\vi,\vj=1,\ldots,\vn$.

Moreover, the derivatives of $\aflow$ are taken at the point $(\vx,\afunc(\vx))$ and
the derivatives of $\bflow$ are taken at $(\aflow(\vx,\afunc(\vx)), \bfunc(\vx))$.
For $\vx=0$ both points coincide with the origin $(0,0)$ since $\afunc(0)=\bfunc(0)=0$ and $\aflow_{0}=\bflow_{0}=\id_{\bnbh}$.
At the origin $(0,0)$ these matrices and vectors have the following form:
$$\pdAdx{\aflow}(0,0)=\pdAdx{\bflow}(0,0)=\mEm, \quad
\pdAdt{\aflow}(0,0)=\aFld(0), \quad \pdAdt{\bflow}(0,0)=\bFld(0).$$ 
Hence 
$$
\frac{d\amap}{d\vx} \ = \
\mEm \ + \ \aFld \cdot \tran{\nabla{\afunc}} 
  \ + \ \bFld \cdot \tran{\nabla{\bfunc}}.
$$
Similar arguments hold for $\vn\geq3$.
Lemma~\ref{lm:Jf0_E_X} and Theorem~\ref{th:f_loc_dif} are proved.
\endproof

\begin{rem}\em
We can imagine that vector fields $\Fld_{\vj}$ define some ``singular coordinate system'' on $\manif$ so that the functions $\afunc_{\vi}$ are ``coordinate functions'' of $\amap$ in this ``coordinate system''.
Thus the matrices $\mEm+\mX$ and $\mEn+\mY$ are the Jacobi matrices of $\amap$ with respect to the usual and ``singular'' coordinate systems respectively. 
We will see in Remark~\ref{rem:Fol_manif} that if this ``coordinate system'' is ``right'', then matrices $\mX$ and $\mY$ coincides.
\end{rem}

\section{Representation of a leaf-preserving mapping as a shift.}\label{sect:repr_sm_shifts}
We give here examples of families of vector fields for which every leaf-preserving diffeomorphism (and even leaf-preserving mapping) can be represented (at least locally) as a smooth shift via some functions.

\subsection{Non-singular foliation.}
Consider the standard non-singular foliation of $\RRR^{\vm}$ whose leaves are $\vn$-dimensional planes defined by the following system of equations: 
$$\vx_{\vn+1}=c_{\vn+1}, \qquad \vx_{\vn+2}=c_{\vn+2}, \qquad \ldots \qquad \vx_{\vm}=c_{\vm},$$
for each $(c_{\vn+1},\ldots,c_{\vm})\in\RRR^{\vm-\vn}$.
Then $\Fol$ is tangent to the following vector fields
$\Fld_{\vi}=\frac{\partial}{\partial\vx_{\vi}}\ (\vi=1,\ldots,\vn)$ that generate the flows $\flow_{\vi}:\RRR^{\vm}\times\RRR\to\RRR^{\vm}$ by:
\begin{equation}\label{equ:nonsing_fol_flow}
\flow_{\vi}(\vx_{1},\ldots,\vx_{\vi}, \ldots, \vx_{\vm},\vt)=
(\vx_{1},\ldots,\vx_{\vi}+\vt, \ldots, \vx_{\vm})
\end{equation}

\begin{lem}\label{lem:non_sing_fol}{\em (cf. Formula (10) of~\cite{Maks:Shifts})}
Let $$\amap=(\amap_1,\ldots,\amap_{\vm}):\RRR^{\vm}\to\RRR^{\vm}$$ be a map preserving each leaf of $\Fol$.
Then $\amap$ is a smooth shift along $\flow_1,\ldots,\flow_{\vn}$ via the functions 
$\afunc_{\vi}=\amap_{\vi}-\vx_{\vi}$ for $\vi=1,\ldots,\vn$.
\end{lem}
\proof 
Since $\amap$ is a leaf-preserving, it follows that 
$\amap_{\vi}(\vx_1,\ldots,\vx_{\vm})=\vx_{\vi}$ if $\vn+1\leq\vi\leq\vm$.
Moreover, if $1\leq\vi\leq\vn$, then
\begin{equation}\label{equ:nonsing_fol_flow_shift_func}
\amap_{\vi}(\vx_1,\ldots,\vx_{\vm})=\vx_{\vi} + (\amap_{\vi}-\vx_{\vi})=\vx_{\vi} + \afunc_{\vi}.
\end{equation}
Hence \
$\amap(\vx) = \flow_{\vn}\bigl(\cdots \flow_{2}(\flow_{1}(\vx,\afunc_1(\vx)), \afunc_{2}(\vx))\cdots,\afunc_{\vn}(\vx) \bigr).$
\endproof

Thus Theorem~\ref{th:f_loc_dif} shows that the Jacobian of $\amap$ can be expressed through the partial derivatives of functions $\amap_{\vi}-\vx_{\vi} \ (\vi=1,\ldots,\vn)$ by the former~$\vn$ coordinates $\vx_{1},\ldots,\vx_{\vn}$.

\begin{rem}\label{rem:Fol_manif}\em
In particular, if $\vm=\vn$, then $\Fol$ consists of a unique leaf $\RRR^{\vm}$ and $\amap$ is an arbitrary diffeomorphism of $\RRR^{\vm}$.
We will see now that in this case the statement of Theorem~\ref{th:f_loc_dif} is trivial.
Indeed,
$$
\dfftotal{\afunc_1,\ldots,\afunc_{\vn}}{\aFld_1,\ldots,\aFld_{\vn}} =
| \mEn + \mY | =
\left|\delta_{ij} + \frac{\partial}{\partial\vx_{\vi}}(\amap_{\vj} - \vx_{\vj}) \right| =
\left| \frac{\partial\amap_{\vj}}{\partial\vx_{\vi}} \right|  =
J(\amap).
$$ 
\end{rem}

\subsection{Foliation generated by product of local flows.}\label{sect:defn_LFP}
Let us represent $\RRR^{\vm}$ as $\RRR^{\vm_1}\times\cdots\times\RRR^{\vm_{\vn}}$, where $\vm=\vm_1+\ldots+\vm_{\vn}$ and each $\vm_{\vi}\geq 0$.
Then every point $\vx\in\RRR^{\vm}$ can be represented in the following form: $\vx=(\vx_{1},\ldots,\vx_{\vm})$, where $\vx_{\vi}\in\RRR^{\vm_{\vi}}$.

For each $\vi=1,\ldots,\vn$ let $\Fld_{\vi}$ be a vector field on $\RRR^{\vm_{\vi}}$. 
We can regard $\Fld_{\vi}$ as a vector field on $\RRR^{\vm}$ that depends only on $\vx_{\vi}$:
$$ 
\Fld_{\vi}(\vx) \equiv (\underbrace{0,\ldots,0}_{\vi-1},\Fld_{\vi}(\vx_{\vi}),0,\ldots,0).
$$
These vector fields define the (in general singular) foliation $\Fol$ on $\RRR^{\vm}$ whose leaves are products of the orbits of $\Fld_{\vi}$.
Thus if $\omega_{\vi}(\vx_{\vi})$ is the orbit of $\Fld_{\vi}$ passing through $\vx_{\vi}\in\RRR^{\vm_{\vi}}$, then 
$$\Fol_{\vx}=\omega_{1}(\vx_1)\times\cdots\times \omega_{\vn}(\vx_{\vn}) \subset \RRR^{\vm}$$
is the leaf of $\Fol$ through $\vx=(\vx_{1},\ldots,\vx_{\vn})$.

Evidently, $\Fol_{\vx}$ is an immersed submanifold of $\RRR^{\vm}$ and its dimension is equal to the total number of those indices $\vi=1,\ldots,\vn$ for which $\Fld_{\vi}(\vx_{\vi})\not=0$.

The following statement is a direct corollary of the results of~\cite[Section~7]{Maks:Shifts}.
\begin{lem}
Suppose that each vector field $\Fld_{\vi}$ is either linear:
$$
\Fld_{\vi}(\vx_1,\ldots,\vx_{\vm_{\vi}}) \ =  \
(\vx_1,\ldots,\vx_{\vm_{\vi}}) \ \matrA_{\vi} \
\left(\begin{matrix}\frac{\partial}{\partial \vx_{1}} \\ \cdots \\ \frac{\partial}{\partial \vx_{\vm_{\vi}}} \end{matrix}\right),
$$
where $\matrA_{\vi}$ is a constant $\vm\times\vm$-matrix, or 
$$
\Fld_{\vi}=\frac{\partial}{\partial \vx_{\vj}}, \qquad \text{where} \qquad  
\vm_{1}+\ldots+\vm_{\vi-1} < \vj \leq \vm_{1}+\ldots+\vm_{\vi},
$$ 
or $\Fld_{\vi}\equiv0$.
Then every smooth leaf-preserving map $$\amap:\RRR^{\vm}\to\RRR^{\vm}$$ is a smooth shift along $\Fld_{1},\ldots,\Fld_{\vn}$ via some smooth functions $\afunc_1,\ldots,\afunc_{\vn}$. 
\end{lem}
\proof
Let $\amap_{\vi}:\RRR^{\vm}\to\RRR^{\vm_{\vi}}$ be the $\vi$-th ``coordinate function'' of $\amap$ and $$\flow_{\vi}:\RRR^{\vm_{\vi}}\times\RRR\to\RRR^{\vm_{\vi}}$$ be the flow generated by $\Fld_{\vi}$.
Since $\amap$ preserves each leaf of $\Fol$ being the products of orbits of $\flow_{\vi}$, it follows that $\amap_{\vi}$ can be regarded as a family of smooth shifts 
\begin{equation}\label{equ:fi_param}
 \amap_{\vi;(\vx_{1},\ldots,\vx_{\vi-1}, -,\vx_{\vi+1},\ldots,\vx_{\vn})}: 
 \RRR^{\vm_{\vi}}\to\RRR^{\vm_{\vi}}
\end{equation}
along the orbits of $\flow_{\vi}$ depending on the parameter that runs over 
$$
\RRR^{\vm_{1}}\times\cdots\times\RRR^{\vm_{\vi-1}}
\times\RRR^{\vm_{\vi+1}}\times\cdots\times\RRR^{\vm_{\vn}}. 
$$ 

Suppose that $\Fld_{\vi}$ is linear and $\matrA_{\vi}\not=0$. 
Then it follows from~\cite[Section~7]{Maks:Shifts} that there exists a smooth function $\afunc_{\vi}:\RRR^{\vm}\to\RRR$ such that 
$\amap_{\vi}(\vx) = \flow_{\vi}(\vx_{\vi}, \afunc_{\vi}(\vx))$.

Suppose that $\Fld_{\vi}=\frac{\partial}{\partial\vx_{\vj}}$.
Then 
$$\flow_{\vi}(\vx_{1},\ldots,\vx_{\vj},\ldots,\vx_{\vm},t) =(\vx_{1},\ldots,\vx_{\vj}+t,\ldots,\vx_{\vm})$$
and 
$\afunc_{\vi}=\amap_{\vj}-\vx_{\vj}$, cf. formulas~\eqref{equ:nonsing_fol_flow} and~\eqref{equ:nonsing_fol_flow_shift_func}, see also~\cite[Formula~(10)]{Maks:Shifts}.

Finally, if $\Fld_{\vi}\equiv0$, then $\flow_{\vi}$ is constant: $\flow_{\vi}(\vx,\vt)=\vx$ and 
$$\amap_{\vj}(\vx)=\vx_{\vj}\qquad \text{for} \qquad \vm_{1}+\ldots+\vm_{\vi-1} < \vj \leq \vm_{1}+\ldots+\vm_{\vi}.$$
Hence we can take $\afunc_{\vi}$ to be an arbitrary smooth function.
\endproof

\end{document}